\input amstex
\documentstyle{amsppt}
\magnification=\magstep1
\define \my{\bf}
\define\po{\parindent 0pt}
\define \mytitle{\po \m \my}
\define \p{{\po \it \m Proof. }}
\define \m{\medskip}
\nologo
\pageheight{16.5cm}
\TagsOnRight
\define \ph#1{\phantom {#1}}
\define \T#1{\widetilde {#1}}

\define \cat{\operatorname{cat}}

\define \eps{\varepsilon}
\topmatter
\title{A remark on fixed point sets of gradient-like flows}\endtitle
\author{Yuli B. Rudyak}\endauthor
\rightheadtext{}
\address{Universit\"at--GH Siegen, FB6/Mathematik,  57068 Siegen,
Germany}
\endaddress 
\email{rudyak\@mathematik.uni-siegen.de}\endemail
\abstract Let $S$ be a set of critical points of a smooth real-valued function on a closed manifold $M$. Generalizing a well-known result of Lusternik--Schnirelmann, Reeken~[R] proved that $\cat S \geq \cat M$. Here we prove a generalization of Reeken"s inequality for gradient-like flows on compact spaces.  
\endabstract
\endtopmatter

\head Introduction\endhead

Given a pair $(X,A)$ of topological spaces, the relative Lusternik--Schnirelmann category of $A$ in $X$, $\cat_XA$, is defined as
$$
\cat_XA:=\min\{k\bigm|A\i U_1\cup \cdots \cup U_{k} \text{ where each $U_i$ is open and contractible in $X$}\}.
$$
Furthermore, we define the Lusternik--Schnirelmann category of a space $X$, $\cat X$, by setting $\cat X:=\cat _XX$, cf. [LS], [F].
\m If $X$ is an absolute neighborhood retract  (ANR) then 
$$
\cat_XA=\min\{k\bigm|A\i F_1\cup \cdots \cup F_{k} \text{ where each $F_i$ is {\it closed} and contractible in $X$}\}.
$$
In particular, $\cat_XA\leq \cat A$ whenever
  $A$ is a closed subset of an ANR $X$. 
\m (The equivalence of this "closed" definition with the "open" one was mentioned by Fox~[F] without proof, so we indicate a proof. Given a closed subset $F$ of an ANR $X$, consider a contracting homotopy 
$$
h_t: F \to F,\quad h_0(f)=f, h_1(f)=x_0 \text{ for every $f\in F$}.
$$
Since $X$ is an ANR. there is an open neighborhood $U$ of $F$ and a homotopy $g_t: U \to X$ which extends $h_t$ and such that $g_0(u)=u$ for every $u\in U$. Finally, we can assume (by taking an open subset of $U$ if necessary) that $g_1(U)$ is contractible in X, because every ANR is locally contractible, [B, IV(3.2)].)
\m Lusternik--Schnirelmann~[LS] proved that, for every smooth function $f: M \to \Bbb R$ on a closed smooth manifold $M$, the number of critical points of $f$ is at least $\cat M$. Reeken~[R] have found the following nice generalization of this result:
$$
\cat S \geq \cat M
$$
where $S$ denotes the set of critical points of $f$.

\m Recently Eugene Lerman asked me whether this result holds in more general situations, for example, if one can prove the similar theorem for orbifolds. In this paper we prove an analog of Reeken"s Theorem for gradient-like flows on compact spaces. 

\m In fact, we consider an axiomatically defined version of Lusternik--Schnirelmann category, and this leads to some new corollaries. For example, we prove that the number of critical levels of any smooth function on a closed manifold $M$ is at least $\cat M$ (this was proved by Takens~[T] under the assumption that all the critical levels of $f$ are path connected).

\head Results \endhead

{\mytitle 1. Definition-Recollection.} A {\it flow} on a topological space $X$ is a family $\Phi=\{\varphi_t\}, t\in \Bbb R$ where each $\varphi_t: X \to X$ is a self-homeomorphism and $\varphi_s\varphi_t=\varphi_{s+t}$ for every $s,t\in \Bbb R$ (notice that this implies $\varphi_0=1_X$). 

\par A flow is called {\it continuous} if the function $X \times \Bbb R \to X, (x,t)\mapsto \varphi_t(x)$ is continuous.
\par A point $x\in X$ is called a {\it rest point} of $\Phi$ if $\varphi_t(x)=x$ for every $t\in \Bbb R$.
\par A continuous flow $\Phi=\{\varphi_t\}$ is called {\it gradient-like} if there exists a continuous function $F: X \to \Bbb R$ with the following property: for every $x\in X$ we have $F(\varphi_t(x))< F(\varphi_s(x))$ whenever $t>s$ and $x$ is not a rest point of $\Phi$. Every such function $F$ is called a {\it Lyapunov function} for $\Phi$.
{\mytitle 2. Definition.} Let $\Phi=\{\varphi_t\}$ be a continuous flow on a topological space $X$. We define a {\it $\Phi$-index function} on $X$ to be any function $\nu: 2^X \to \Bbb N\cup\{0\}$ with the following properties:
\roster \item (monotonicity) If $A\i B \i X$ then $\nu (A)\leq \nu (B)$;
\item (continuity) For every $A\i X$ there exists an open neighborhood $U$ of $A$ such that $\nu (A) =\nu (U)$;
\item (subadditivity) $\nu(A\cup B) \leq \nu (A) + \nu (B)$;
\item (invariance) $\nu (\varphi_t(A))=\nu(A)$ for every $A\i X$ and $t\in \Bbb R$.
\endroster
We say that an index function $\nu$ satisfies the {\it normalization condition} if the following condition holds:
\roster
\item[5] (normalization) $\nu(\emptyset)=0$, $\nu(\{x\})=1$ for every point $x\in X$. 
\endroster
\proclaim{3. Theorem} Let $\Phi$ be a gradient-like flow on a compact metric space $X$, let $S$ be the set of rest points of $\Phi$, and let $F$ be a Lyapunov function for $\Phi$. Suppose that the set $F(S)$ is finite, say, $F(S)=\{a_1, \dots, a_m\}\i \Bbb R$. Then
$$
\sum_{i=1}^m\nu(S\cap F^{-1}(a_i))\geq \nu(X).
$$
Moreover, if $\nu$ satisfies the normalization condition then $m\geq \nu(X)$.
\endproclaim

\p We let $F^a:=F^{-1}(- \infty, a]$. We set
$$
c_k:=\inf \{ c \in \Bbb R \bigm| \nu(F^c)\ge k\},\quad k=0,\ldots, \nu(X).
$$
We notice that that
$$
\nu(F^{c_k-\eps})<k \text{ for every $\eps >0$}, \quad \nu(F^{c_k})\geq k.
$$ 
The first inequality follows from the definition of $c_k$. 
Furthermore, by monotonicity, $\nu (F^{c_k+\eps}) \ge k$ for every $\eps >0$. Moreover, because of compactness of sets $F^d$, for every neighborhood $U$ of $F^{c_k}$ there exist $\eps>0$ such that $F^{c_k+\eps}\subset U$. Thus, by continuity, $\nu (F^{c_k}) \ge k$.    

In particular, every $c_k$ is a critical level of $F$, i.e., $c_k\in F(S)$. Indeed, otherwise we can deform $F^{c_k}$ to  $F^{c_k-\eps}$, and hence $\nu(F^{c_k})\leq \nu(F^{c_k-\eps})$. This contradicts to what we said above.     

\proclaim{4. Lemma{\rm (cf. [LS], [DNF])}} If $c_k=c_{k+1}=\cdots = c_{k+r}$, then $\nu(S\cap F^{-1}(c_k))\geq r+1$. 
\endproclaim
\p Let $A:=S\cap F^{-1}(c_k)$. Then, because of continuity, there is a neighbourhood $U$ of $A$ with $\nu(U)=\nu(A)$. Furthermore, there are $\eps>0 $ and $t\in \Bbb R$ such that $\varphi_t(F^{c_k+\eps}\setminus U)\i F^{c_k-\eps}$. So, $\nu(F^{c_k+\eps}\setminus U)<k$. Now, because of subadditivity,
$$
\nu(F^{c_k+\eps}\setminus U)+\nu(U)=\nu(F^{c_{k+r}+\eps}\setminus U)+\nu(U)\geq \nu(F^{c_{k+r}+\eps})\geq k+r
$$
and so $\nu(U)>r$. Thus, $\nu(A)\geq r+1$. The Lemma is proved.

\m We continue the  proof of the Theorem. We have
$$
c_1=\cdots=c_{i_1}<c_{i_1+1}=\cdots=c_{i_2}<\ldots <c_{i_n+1}=\cdots=c_{i_{n+1}}=c_{\nu(X)}.
$$ 
Now,  set $S_k:=S\cap F^{-1}(c_{i_k})$. Then, by the Lemma, $\nu(S_k)\geq i_k-i_{k-1}$ (where $i_0=0)$. Thus
$$
\sum_{k=1}^{n+1}\nu(S_k)\geq \nu(X).
$$
Now the result follows, because $c_i\in \{a_1,\ldots, a_m\}$ for every $i$.
\m Finally, if $\nu$ satisfied the normalization condition then, by the Lemma and subadditivity, the number of rest points belonging to $S_k$ is at least $i_k-i_{k-1}$, and the result follows.
\qed

\proclaim{5. Corollary}Let $f$ be a smooth function on a closed manifold $M$, let $\{a_1, \ldots,a_m\}$ be the $($finite$)$ set of the critical levels of $F$, and let $S$ be  the set of the critical points of $f$. Then
$$
\sum_{k=1}^m \cat_M(S\cap f^{-1}(a_k) \geq \cat M.
$$
\endproclaim 

\p Given $A\in M$, set $\nu(A)=\cat_M A$ where each $A_i$ is open and contractible in $M$. We leave it to the reader to check that $\nu$ is an index function. Clearly, $\nu(M)=\cat M$.
\qed

\m Since $\cat_M F\leq \cat F$ for every closed subset $F$ of $M$, Corollary 5 implies the inequality 
$$
\cat S\geq \cat M
$$ 
of Reeken [R]. 

\proclaim{Corollary}Let $f$ me a smooth function with isolated critical points on a closed connected manifold $M$. Then the number of critical levels of $f$ is at least $\cat M$.
\endproclaim
\p We use the notation from Theorem 4. Since $M$ is connected, every set $S\cap f^{-1}(a_k)$ is contractible in $M$, and so $\cat_M(S\cap f^{-1}(a_k)=1$. Thus. $m\geq \cat M$.
\qed.

\m If all the critical levels of $f$ are path connected, this results is actually contained in the paper of Takens~[T, Prop. 2.8].


\Refs\nofrills{References}

\widestnumber\key{DNF}

\ref\key{B}
\by K. Borsuk
\book Theory of Retracts
\publ PWN -- Polish Scientific Publischers, Warszawa
\yr 1967
\endref\vskip 6pt

\ref\key{DNF}
\by B. A. Dubrovin, S. P. Novikov, A.T. Fomenko
\book Modern Geometry
\publ  Springer, Berlin Heidelberg New York
\yr 1985
\endref\vskip 6pt

\ref\key{F}
\by R. Fox 
\paper On the Lusternik--Schnirelmann category
\jour Ann. of Math.
\vol 42
\pages 333--370
\yr 1941
\endref\vskip 6pt

\ref\key{LS}
\by L. A. Lusternik, L. G. Schnirelmann
\book Methodes  topologiques dans le probl\`emes variationels
\publ Hermann, Paris 
\yr 1934
\endref\vskip 6pt

\ref\key{R}
\by M. Reeken
\paper Stability of critical points under small perturbations. Part I. Topological theory
\jour Manuscripta Math.
\vol 7
\pages 387--411
\yr 1972
\endref\vskip 6pt

\ref\key{T}
\by F. Takens
\paper The minimal number of crirical points of a function on a compact manifold and the Lusternik--Schnirelmann category
\jour Invent. Math
\vol 6
\yr 1968
\pages 199--224
\endref
\vskip6pt

\endRefs
\enddocument